\documentclass{amsart}

\title{On some properties of Perron numbers}

\author{Nikita Sidorov}
\address{I am retired and have no affiliation. Please use my email for communication.}
\email{nikita.a.sidorov\@ gmail.com}
\dedicatory{To the memory of Vladimir Komarov, one of the bravest men in
human history, who sacrificed his own life to save the life of his friend Yuri Gagarin, thus 
becoming the first human who died in space.}
\date{November 2, 2023}

\begin{document}

\begin{abstract}
Let $\theta$ be a real number, $n\in\mathbb N$, and
\[
D_n(\theta)=\left\{\sum_{k=1}^n a_k\theta^{k}\mid a_k\in\{0,\dots,\lfloor \theta\rfloor\}\right\}.
\]
Let $\theta$ be a Perron number, that is,
an algebraic integer $>1$ whose other Galois conjugates are less than $\theta$
in absolute value. In this note I prove three results. 
\begin{enumerate}
\item\[
\theta^n\ll\#D_n(\theta)\ll\sqrt n \theta^n.
\]
\item If all Galois conjugates of $\theta$ are real, then $\#D_n(\theta)\asymp\theta^n$.
\item $\theta$ is of height $\le \lfloor\theta\rfloor$, that is, $\theta$ is a root of a polynomial
with integer coefficients bounded by $\lfloor\theta\rfloor$ in absolute value.
\end{enumerate}
\end{abstract}

\maketitle

Let $\theta\in\mathbb R$. Put
\[
D_n(\theta)=\left\{\sum_{k=1}^n a_k\theta^{k}\mid a_k\in\{0,\dots,\lfloor \theta\rfloor\} \right\}.
\]
Assume henceforth that $\theta$ is {\it Perron}, that is,
an algebraic integer $>1$ whose other Galois conjugates are less than $\theta$
in absolute value.

\medskip\noindent
{\bf Theorem 1.} {\sl We have}
\[
\theta^n\ll\#D_n(\theta)\ll n^{1/2} \theta^n.
\]
Furthermore, if all Galois conjugates of $\theta$ are real, then $\#D_n(\theta)\asymp\theta^n$.

\medskip\noindent
{\bf Theorem 2.} {\sl $\theta$ is of height $\le \lfloor\theta\rfloor$, that is, $\theta$ is a root of a polynomial
with integer coefficients bounded by $\lfloor\theta\rfloor$ in absolute value.}

\medskip\noindent
I shall prove both statements simultaneously.

\medskip\noindent
  Let $x=\sum_{1}^{n}a_k\theta^k, y=\sum_{1}^k b_k\theta^k$ with
$a_k\in\{0,1,\dots,\lfloor\theta\rfloor\},\ 1\le k\le n$. 
Put $c_k=a_k-b_k\in\{-\lfloor\theta\rfloor,\dots,\lfloor\theta\rfloor\}$. Define
\[
S_n(\theta)=\sum_{k=1}^{n}c_k\theta^k.
\]  

Put $\alpha_k(\theta)=\sum_{j=1}^{d}\theta_j^k\in\mathbb N\sim\theta^k$.
Assume first that $\theta$ has a non-real Galois conjugate.
Denote the non-real Galois conjugates of $\theta$ 
by $\theta_2,\dots, \theta_s$ and the real ones by $\theta_{s+1},\dots, \theta_d$. 
Recall that if $\theta_j\in\mathbb{C}\setminus \mathbb{R}$, then $\overline{\theta_j}$ is a Galois conjugate of $\theta$ as well.  
(If all the Galois conjugates of $\theta$ are real, just ignore the complex part
of the argument.)

Since all non-real $\theta_j$ are irrational, by Birkhoff's ergodic theorem which
in the case of a uniquely ergodic irrational rotation by $\theta_j/|\theta_j|$ applies to {\bf all} arguments,
\[
\lim_{n\to\infty} \frac{1}{n}\sum_{k=1}^{n}\text{Re}(\theta_j^k/|\theta_j^k|)=\int_{-1}^{1} x\ dx=0,
\]
whence by the Central Limit Theorem [3],
\[
\sum_{k=1}^{n} \text{Re}(\theta_j^k)= O\left(n^{1/2}|\theta_j|^n\right),
\]
Therefore,
\[
\sum_{k=1}^{n}\sum_{j=1}^{s}\text{Re}(\theta_j^k)=
O\left(n^{1/2}\max_{j\in\{1,\dots,s\}}|\theta_j^n|\right),
\]
and 
\[
\sum_{k=1}^{n}\sum_{j=s+1}^{d} \theta_j^k = O\left(\max_{j\in\{s+1,\dots,d\}}{|\theta_j^n|}\right).
\]
Summing up,
\[
\sum_{k=1}^{n}\sum_{j=1}^{d} \theta_j^k = O\left(n^{1/2}\max_{j\in\{1,\dots,d\}}{|\theta_j^n|}\right).
\]
Theorem 1 is now proved since $\#D_n(\theta)\gg \theta^n$ for all $\theta>1$ -- just count the prefixes of length~$n$ for the greedy $\theta$-expansions
[4]. If all Galois conjugates of $\theta$ are real, it is clear from the above argument that $\#D_n(\theta)\ll\theta^n$.

The claim of Theorem~2 is a trivial consequence of that of Theorem~1.
Indeed, if the height of $\theta$ is $\ge\lfloor\theta\rfloor+1$,
then $\#D_n(\theta)=(\lfloor\theta\rfloor+1)^n$, which contradicts
Theorem~1.

\medskip\noindent
\textbf{Remarks. 1.}
Thus, the Perron case is very similar to the Pisot case (see, {\it e.g.} [2]), with an extra $\sqrt n$ multiplier thrown in the mix -- which does not affect the corresponding arguments.

\medskip\noindent\textbf{2.} As $n\to\infty$, the minimum of $|x-y|$ with
the same notation as above, is known to be $o(\theta^{-n})$ -- see [1].


\begin{thebibliography}{9}

\bibitem{[1]} D.-J. Feng, {\it On the topology of polynomials with bounded integer coefficients},  J. Eur. Math. Soc. {\bf 18} (2016), 181--193.

     
\bibitem{[2]} K. G. Hare, Ph.~D. Thesis, https://uwaterloo.ca/scholar/kghare/publications/pisot-numbers-and-spectra-real-numbers
 
\bibitem{[3]} https://mathoverflow.net/questions/172936/central-limit-theorems-for-irrational-rotation

\bibitem{[4]} A. R\'enyi, {\em Representations for real numbers and
their ergodic properties,} Acta Math. Acad. Sci. Hung. {\bf 8}
(1957) 477--493.


\end{thebibliography}
\end{document}